\documentclass[twocolumn,10pt]{extarticle}
\usepackage{abstract}
\usepackage[default]{raleway}
\usepackage[T1]{fontenc}

\usepackage{hyperref}

\usepackage[margin=0.875in]{geometry}
\usepackage{placeins}

\usepackage[parfill]{parskip}
\usepackage{amsmath} 
\usepackage[utf8]{inputenc}

\usepackage{graphicx}
\usepackage{multirow}

\usepackage[numbers]{natbib}

\begin{document}

\title{Softplus Penalty Functions for Constrained Optimization}
\author{
  Stefan Meili\thanks{Stefan Meili \href{mailto:smeili@noram-eng.com}{smeili@noram-eng.com}, NORAM Engineering and Constructors Ltd., Vancouver, Canada}
}
\date{}

\twocolumn[
\maketitle
\begin{onecolabstract}
\vspace{12pt}
Penalty functions are widely used to enforce constraints in optimization problems and reinforcement leaning algorithms. Softplus and algebraic penalty functions are proposed to overcome the sensitivity of the Courant-Beltrami method to strong objective function gradients. These penalty functions are shown to converge in one fifth as many iterations on problems with fifty optimized parameters and be less sensitive to objective function gradient and scaling than the Courant-Beltrami penalty function. The softplus and algebraic penalty functions developed in this work allow optimized solutions to be found in fewer iterations with greater control over constraint accuracy.
\vspace{24pt}
\end{onecolabstract}
]

\saythanks

\section{Introduction}

Optimization and reinforcement learning problems problems are found throughout science, engineering, machine learning, economics, and operations research. This class of problems attempts to find the set of variables, $U$, that produces the best solution possible. Solution quality is judged using an objective function, $O(U)$, where a minimum value indicates the best solution.

Penalty functions are a simple method of enforcing constraints on an optimization problem. Infeasible solutions to the objective function are penalized by returning a large value when a constraint is violated. In this way, a constrained optimization problem can be solved using an unconstrained optimization algorithm.

To implement, constraint error $x$ is calculated using Equation \ref{eq:const_error}, where $v$ is a quantity derived from the solution space $U$, and passed to a penalty function such as Equation \ref{eq:penalty_lin} or \ref{eq:penalty_cb}. 
\begin{align}\label{eq:const_error}
x = v - v_{\text{target}}
\end{align}
\begin{align}\label{eq:penalty_lin}
g^{+}_{\text{Lin}}(x) = max(0,x)
\end{align}
\begin{align}\label{eq:penalty_cb}
g^{+}_{\text{C-B}}(x) = max(0,x)^2
\end{align}
\begin{align}\label{eq:psum}
P_{\text{sum}}(U) = \sum_{i=0}^{n} \sigma_i \cdot g_i(x_i(U))
\end{align}
\begin{align}\label{eq:objective}
O_c(U) = O(U) + P(U)
\end{align}
Multiple penalty functions are summed using Equation \ref{eq:psum} and are often weighted with scale parameter $\sigma$ to improve adherence to each constraint \cite{evolutionary_comp}. The constrained objective function, Equation \ref{eq:objective}, can then be passed to the optimization algorithm.

Applying a linear penalty function, Equation \ref{eq:penalty_lin}, creates a sharp corner as the constraint becomes active which cannot be differentiated. This can stall or destabilize convergence of gradient descent algorithms.

The Courant-Beltrami penalty function, Equation \ref{eq:penalty_cb}, is perhaps the most commonly used penalty function. It has a derivative of zero at $x=0$, allowing infeasible solutions when applied to an objective function with a strong gradient.  Increasing $\sigma$ reduces constraint violations, but can lead to machine overflow in practice. Furthermore, this approach becomes unstable when used with gradient based optimization algorithms as $1/\sigma$ approaches the interval, $\Delta U$, used to approximate the local gradient of the objective function.

Smooth approximations of $max(x, 0)$ such as the one-sided huber loss function can reduce the computational complexity of constrained optimization problems with many dimensions \cite{tatarenko2020smooth} \cite{nedich2020convergence}. However, this penalty function is locally quadratic near constraints and relies on scaling to enforce constraints in the presence of strong gradients.

This work proposes variations of the softplus function \cite{benigo} for use as penalty functions. These functions are continuously differentiable, have a slope at $x=0$, and can produce solutions with small constraint errors with low computational effort.

\FloatBarrier
\section{Softplus Penalty}\label{softplus}

The softplus function is derived by integrating the logistic function, $\frac{1}{1+e^-x}$. Other sigmoid functions such as the algebraic expression $\frac{x}{\sqrt{1 + x^2}}$ produce similar approximations.

For this work, parameters $\beta^{-}$ and $\beta^{+}$ were introduced to set the negative and positive saturation values respectively. A hardness parameter, $\alpha$, controls how quickly the sigmoid functions transition between the negative and positive saturation values. See Equations \ref{eq:siglog} and \ref{eq:sigalg}.

\begin{align}\label{eq:siglog}
Logistic(x) = \frac{\beta^{+}-\beta^{-}}{1 + 2^{\frac{-x}{\alpha}}} +\beta^{-}
\end{align}

\begin{align}\label{eq:sigalg}
Algebraic(x) = \frac{x \cdot (\beta^{+}-\beta^{-})}{2 \sqrt{4 \alpha^2 +  x^2 }} + \frac{\beta^{-}+\beta^{+}}{2}
\end{align}

Table \ref{tab:penalty} lists penalty functions derived from Equations  \ref{eq:siglog} and \ref{eq:sigalg} by integration. To generate a penalty function that approximates $max(0,x)$ and enforces a less-than constraint, $\beta^{-}=0$ and $\beta^{+}=1$. For an equality, $\beta^{-}=-1$ and $\beta^{+}=1$, while for a greater-than constraint, $\beta^{-}=-1$ and $\beta^{+}=0$.

\begin{table}[h]
\caption{Smooth Penalty Functions}\label{tab:penalty}
\rule{0pt}{0.5ex}
\centering
\resizebox{\linewidth}{!}{
\begin{tabular}{ c | c c }
  \rule{0pt}{2ex} $g(x, \alpha)$ & Algebraic & Softplus \\ [1ex]
  \hline
  \rule{0pt}{4ex} $g^{+}, g^{<}$ & $\displaystyle\frac{\sqrt{4 \alpha^2 + x^2} + x }{2}$ & \
  $\displaystyle\alpha \log_{2}(1+2^\frac{x}{\alpha})$ \\[2ex]
 
  \rule{0pt}{2ex} $g^{0}, g^{=}$ & $\displaystyle\sqrt{4 \alpha^2 + x^2}$ & \
  $\displaystyle 2 \alpha \log_{2}(1+2^\frac{x}{\alpha}) - x$ \\[2ex]

  \rule{0pt}{2ex} $g^{-}, g^{>}$ & $\displaystyle\frac{\sqrt{4 \alpha^2 + x^2} - x }{2}$ & \
  $\displaystyle\alpha \log_{2}(1+2^\frac{-x}{\alpha})$\\[1ex]
\end{tabular}}
\end{table}

Combining linear constraints using Equation \ref{eq:psum} introduces a mitred effect in regions where multiple constraints overlap. Combining constraints using the euclidean norm, Equation \ref{eq:pnorm}, was found to improve convergence, as it offers a more gradual transition between adjacent constraints.

\begin{align}\label{eq:pnorm}
P_{norm}(U) = \sqrt{\sum_{i=0}^{n} \left[\sigma_i \cdot g_i(x_i(U),\alpha_i)\right]^2}
\end{align}

\begin{figure}[h]
  \includegraphics[width=\linewidth]{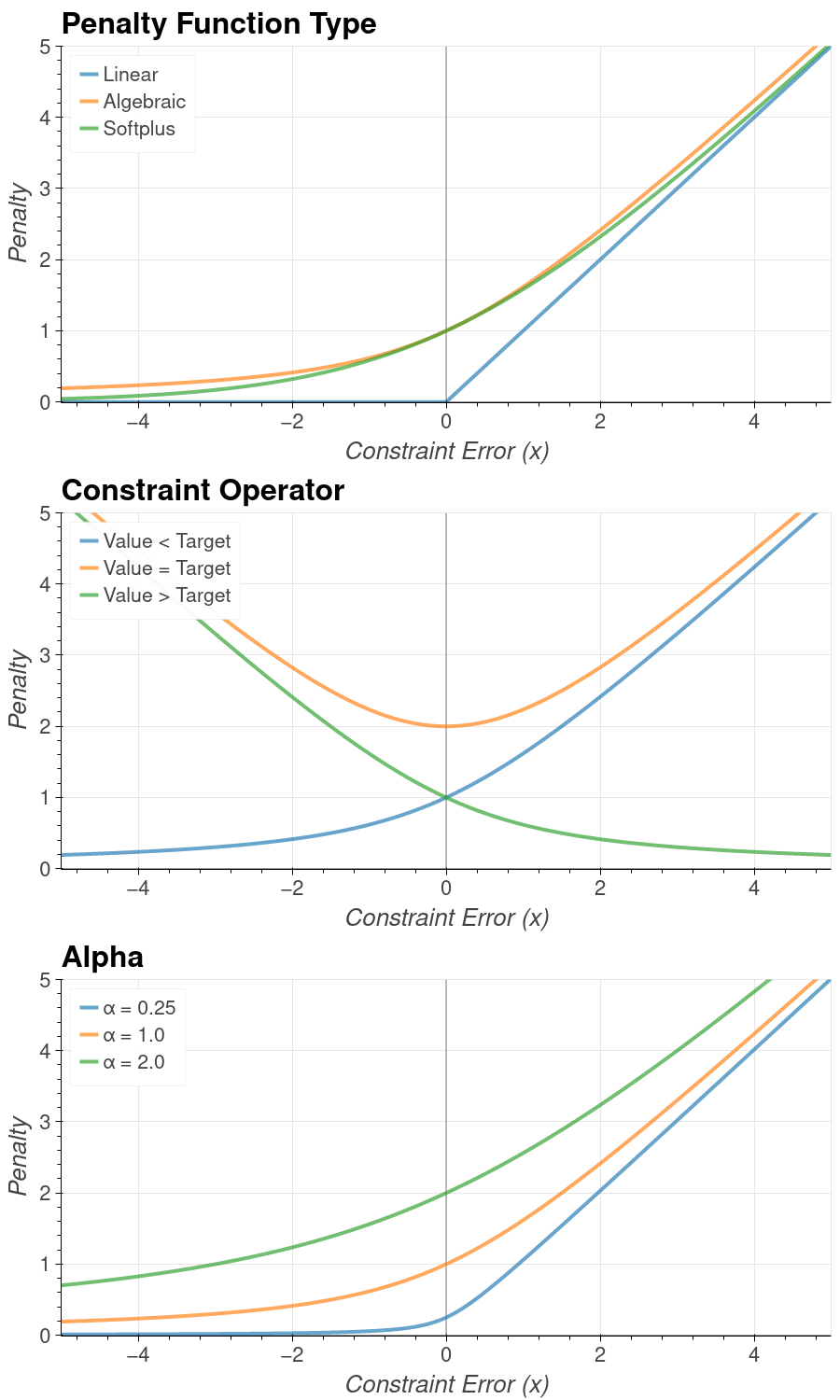}
  \caption{Penalty Function Parameters}
  \label{fig:parameters}
\end{figure}

Figure \ref{fig:parameters} compares the algebraic and softplus penalty functions to $max(0,x)$. Softplus functions are characterized by a round corner near $x=0$ and quick decay to $max(0,x)$. The algebraic functions have a harder corner and slow decay. Corner hardness in the region of $x=0$ is controlled by $\alpha$. Small values improve accuracy, but overflow and optimization convergence can become problematic as the corner sharpens.

In practice, the softplus penalty functions in Table \ref{tab:penalty} are approximately 5X slower than the algebraic functions when implemented in Python and Numpy. Furthermore, the  $2^\frac{x}{\alpha}$ term in the softplus functions overflows at small values. For 64bit, $x/\alpha < 1024$, and for 32bit, $x/\alpha < 128$. As $2^{\frac{x}{\alpha}}$ approaches overflow, the softplus penalty functions should be masked to return an appropriate linear response.

\FloatBarrier
\section{Parameters}\label{parameters}

\begin{figure}[hbt!]
  \includegraphics[width=\linewidth]{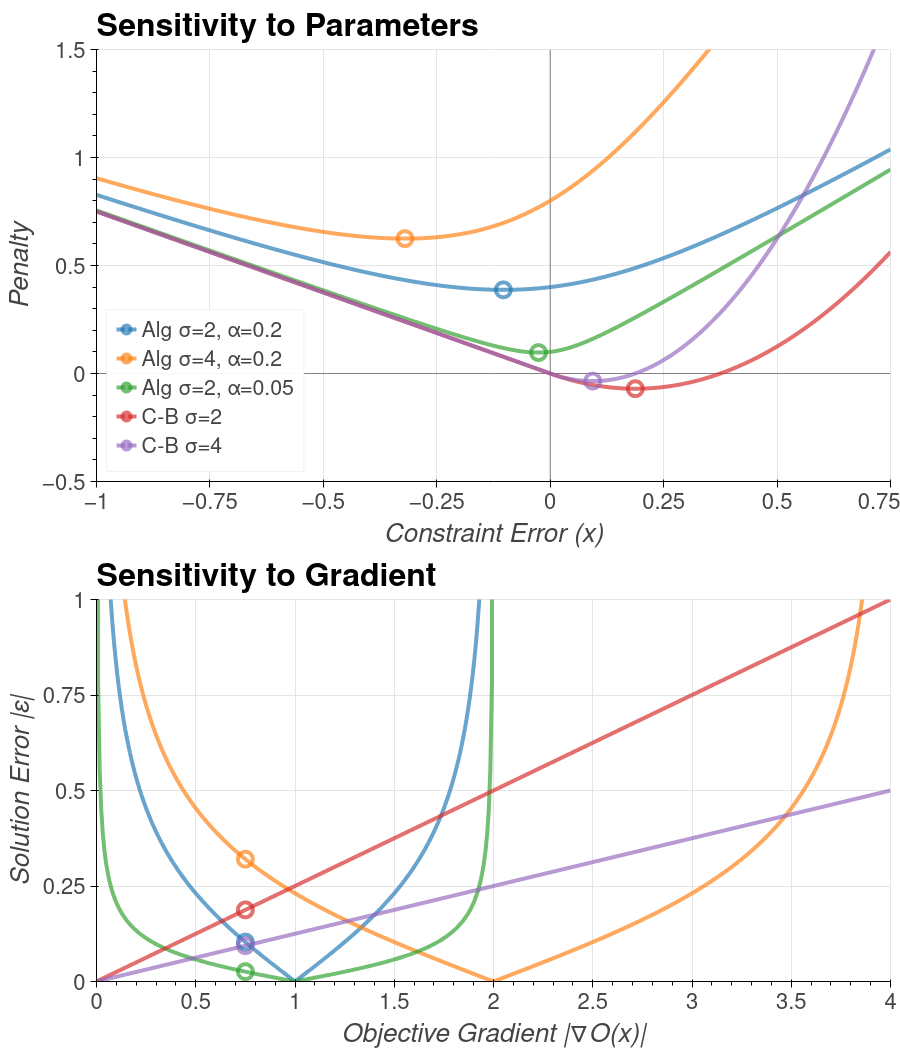}
  \caption{Sensitivity to Parameters and Gradient}
  \label{fig:accuracy}
\end{figure}

Figure \ref{fig:accuracy} explores the effect of parameters $\alpha$ and $\sigma$ and objective function gradient, $\left|\nabla O(x)\right|$, on the solution error, $|\epsilon|$, of the algebraic and Courant-Beltrami penalty functions. Increasing $\sigma$ increases solution error of the algebraic penalty function, and reduces the error of the Courant-Beltrami penalty function. Reducing $\alpha$ reduces the error produced when using the algebraic and softplus penalty function, as the minimum is drawn closer to the true constraint boundary.

The error produced when using the Courant-Beltrami error functions increases proportionally with objective function gradient. The algebraic penalty function produces zero solution error when $\sigma = 2|\nabla O(x)|$ and increases asymptotically as $x$ approaches 0 and 1.

\begin{figure}[hbt!]
  \vspace{36pt}
  \includegraphics[width=\linewidth]{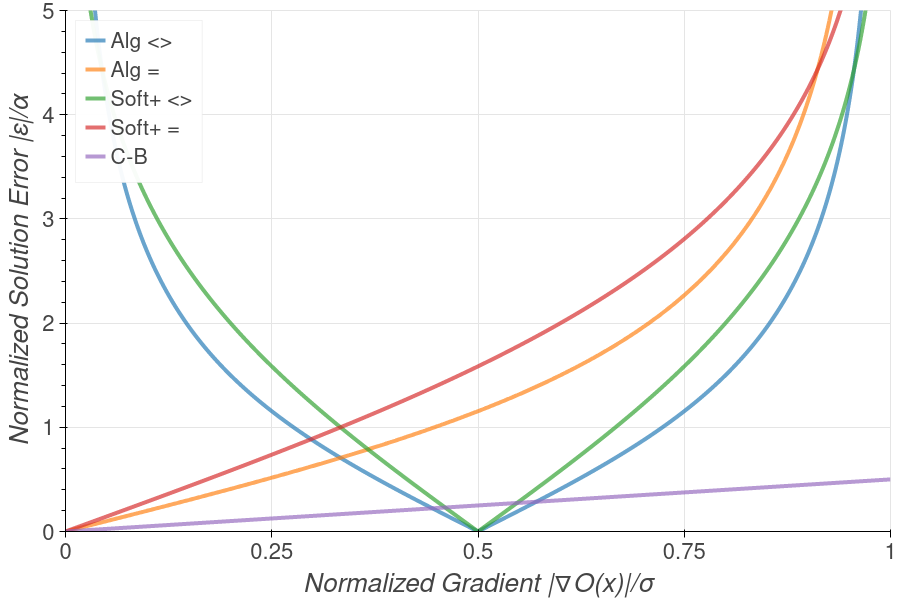}
  \caption{Normalized Penalty Function Error}
  \label{fig:norm_err}
\end{figure}

Solution error can be derived using Equation \ref{eq:err_deriv} if the local gradient is known near the optimized solution. Equations \ref{eq:err_cb} - \ref{eq:err_log=} estimate solution error of the Courant-Beltrami, algebraic and softplus penalty functions and are plotted in Fig. \ref{fig:norm_err}.

\begin{align}\label{eq:err_deriv}
\frac{\partial}{\partial x} \left( P(x) - \nabla x \right )= 0
\end{align}
\begin{align}\label{eq:err_cb}
\left| \epsilon_{\text{C-B}} \right| = \left|\frac{\nabla}{2 \sigma}\right|
\end{align}
\begin{align}\label{eq:err_alg<>}
\left| \epsilon_{\text{Algebraic <>} } \right| =\left| \alpha\sqrt{\frac{1}{\nabla (\sigma-\nabla)}} (\sigma - 2 \nabla) \right|
\end{align}
\begin{align}\label{eq:err_alg=}
\left| \epsilon_{\text{Algebraic =}} \right| = \left| 2 \nabla \alpha \sqrt{\frac{1}{(\sigma-\nabla) (\nabla+\sigma)}}\right|
\end{align}
\begin{align}\label{eq:err_log<>}
\left| \epsilon_{\text{Softplus <>}} \right| = \left| \alpha \log_2\left(\frac{\nabla}{\sigma-\nabla}\right) \right|
\end{align}
\begin{align}\label{eq:err_log=}
\left| \epsilon_{\text{Softplus =}} \right| = \left| \alpha \log_2\left(\frac{\nabla + \sigma}{\sigma-\nabla}\right) \right|
\end{align}

The Courant-Beltrami, algebraic equality and softplus equality penalty functions converge to zero solution error only when the objective function gradient is zero. The algebraic and softplus inequality penalty functions produce zero solution error when $\sigma$ is twice the local gradient. This suggests an adaptive algorithm where $\sigma$ is set to twice the gradient between iterations.

The algebraic functions perform poorly in comparison to the softplus functions at extreme values of $|\nabla O(x)/\sigma|$ near zero and one due to their slow approach to the linear ideal. The algebraic functions have a relatively strong slope away from $x=0$ that can dominate a shallow gradient. When $|\nabla O(x)| > \sigma$, the softplus and algebraic penalty functions are ineffective. This limitation doesn't apply to the Courant-Beltrami penalty function, as the response continues to grow as $x\rightarrow \infty$.

Softplus, algebraic  and Courant-Beltrami penalty functions produce similar average solution error over the range of 0.05 $< |\nabla O(x)/\sigma| <$ 0.95 when $\alpha \approx$ 0.1. Equations \ref{eq:err_alg<>} and \ref{eq:err_log<>} have the same slope near $|\nabla O(x)/\sigma| =$  0.5 when $\alpha = log(2)$ for the softplus function.

\begin{figure}
  \includegraphics[width=\linewidth]{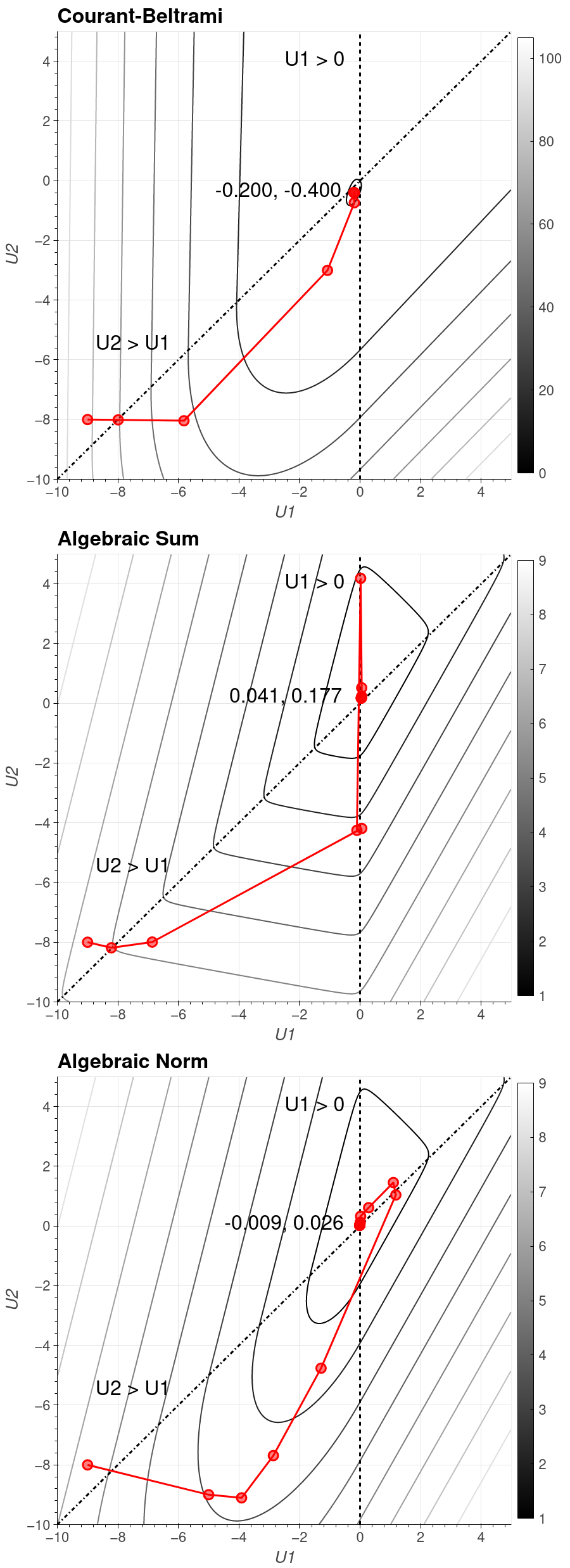}
  \caption{Convergence Behaviour}
  \label{fig:convergence}
\end{figure}

Figure \ref{fig:convergence} compares the convergence behaviour of various penalty functions when passed to the BFGS optimization algorithm \cite{practical_opt} implemented in Scipy 1.6.2. $O(U)=0.2 U_1 + 0.2 U_2$ is optimized subject to $U_1>0$ and $U_2 > U_1$ with minimum at $[0,0]$.

The Courant-Beltrami penalty function converges quickly, but with significant error due to the sensitivity of this penalty function to solution gradient. The optimization paths produced by using the algebraic penalty functions with $\alpha =$ 0.1 are also illustrated. When summed using Equation \ref{eq:psum}, the optimization path overshoots constraints, but reliably converges to a small error. Applying  the euclidean norm, Equation \ref{eq:pnorm}, yields steady convergence to very low error.

Figure \ref{fig:convergence} illustrates a disadvantage of combining penalty functions with the euclidean norm where constraints intersect at acute angles. In these regions, the local gradient of the penalty function becomes very shallow and may not provide an effective constraint against the objective function. Increasing $\sigma$, summing constraints with \ref{eq:psum}, or adding the Courant-Beltrami penalty function can help.

Parameters for the softplus and algebraic penalty functions found in Table \ref{tab:penalty} should be selected as follows
\begin{itemize}
\item $\sigma$ should be set a small multiple of the gradient if it is approximately known to provide some margin before the constraints become ineffective if the gradient is underestimated.
\item If the local objective function gradient is known, setting $\sigma = 2|\nabla O(U)|$ produces solutions with very low error.
\item $\alpha$ directly controls solution error and can approach the jacobian approximation interval before destabilizing optimization.
\item Combining softplus and algebraic penalty functions by summation or adding Equation \ref{eq:penalty_cb} can prevent them from becoming ineffective when constraints intersect at acute angles.
\item The Algebraic functions perform poorly in problems where a gentle gradient is present or when scaled very aggressively. Softplus functions are recommended in these applications.
\end{itemize}

\section{Experiments}\label{experiments}

The performance of the penalty functions developed in Section \ref{softplus} were measured on two optimization problems with linear objectives.

In the first problem, sheared hyperplane constraints were applied to investigate the performance of the penalty functions in corners. For $N$ dimensions, $2N$ constraints were applied at a positive and negative coordinate of each dimension to form a hypercube. The resulting planes were then sheared by a random factor $int(N/2)$ times. This problem is guaranteed to converge in a corner with $N$ active constraints, and the exact solution can be found analytically.

In the second problem, a single hyperspherical constraint was applied to simulate convergence against a single curved constraint. The exact solution is trivial to determine analytically.

The objective gradient was randomly selected between 1e-2 and 5 and an initial point was randomly selected in each experiment. The BFGS optimization algorithm in Scipy 1.6.2 was used to find the optimal solution. The Jacobian was estimated using a central second order approximation with a fixed interval size of 1e-6. 500 samples were run for each problem type, dimension and penalty function configuration. For each experiment, final solution error $|\epsilon| = |U_{opt}-U_{true}|$ and optimization algorithm iterations were recorded.

Four penalty function configurations were explored: algebraic penalty functions were combined by summation and the euclidean norm, while the Courant-Beltrami penalty functions were summed, and the softplus penalty functions were combined with the norm. Penalty function parameters were selected to produce similar solution error on the sheared hyperplane problem. $\sigma =$ 15 and $\alpha =$ 3e-5 was set for the softplus and algebraic experiments, and $\sigma =$ 1e4 for the Courant-Beltrami experiments.

A total of four sets of experiments were run. In the first two sets, the Courant-Beltrami, algebraic, and softplus penalty functions were compared on the sheared hyperplane and hypersphere problems with 2 to 50 dimensions. In a further two experiments, the effect of $\sigma$ was assessed at 12 dimensions between 10 and 1e7, and the effect of $\alpha$ assessed between 1e-7 and 0.1.

The experiment source code is hosted at: \href{https://github.com/stefanmeili/softplus-penalty-functions}{https://github.com/stefanmeili/softplus-penalty-functions}.

\section{Results}\label{results}

\begin{figure}
  \includegraphics[width=\linewidth]{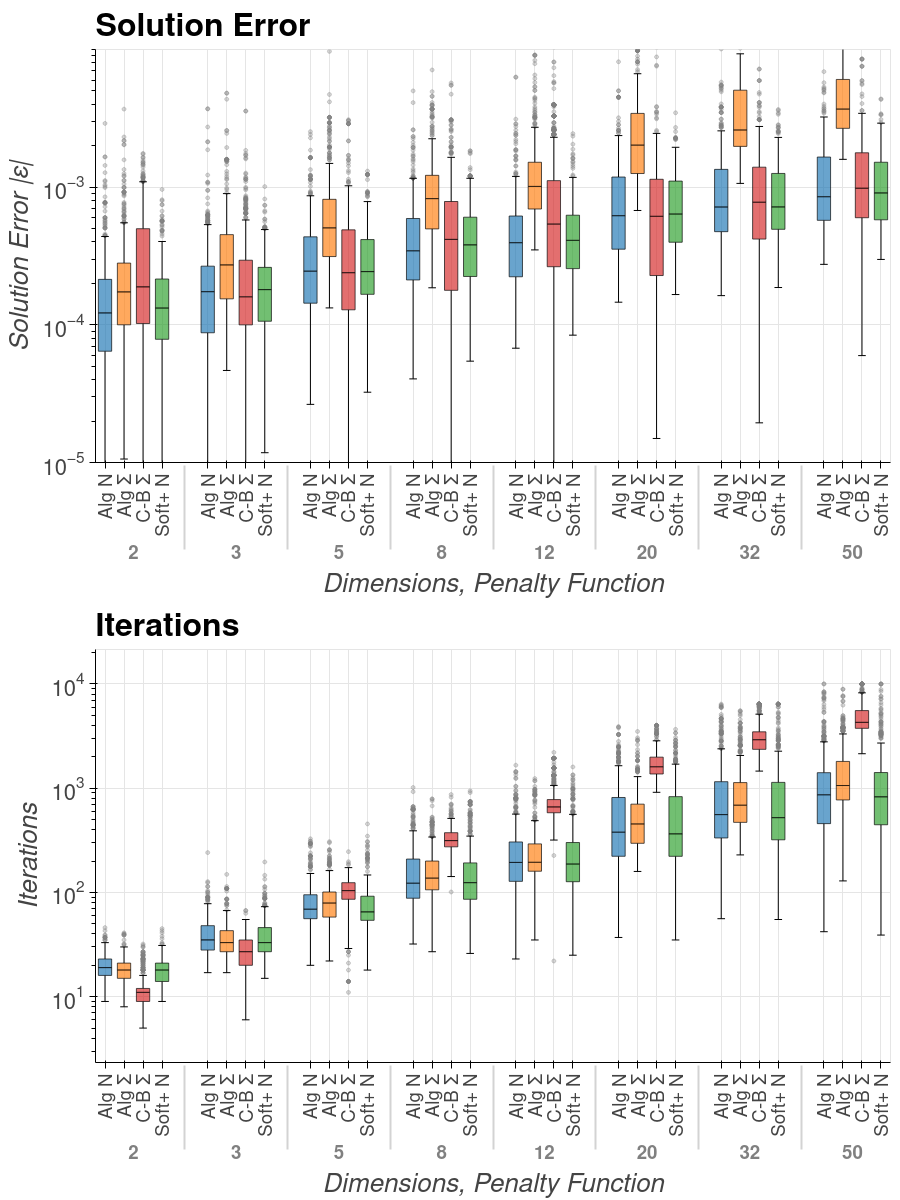}
  \caption{Results - Sheared Hyperplanes}
  \label{fig:exp_dims_shear}
\end{figure}

Figure \ref{fig:exp_dims_shear} and Table \ref{tab:exp} summarize the outcome of the first experiment on the sheared hyperplanes problem. Algebraic and softplus penalty functions combined with the norm and Courant-Beltrami penalty functions converge to solutions with comparable constraint error ranging from approximately 1.5e-4 at two dimensions to 9e-4 at 50 dimensions. Summing algebraic penalty functions produced errors 4 times larger at 50 dimensions.

The results of the sheared hyperplane experiments demonstrate that the algebraic and softplus penalty functions offer a clear advantage in computational efficiency. At 50 dimensions, the Courant-Beltrami penalty functions required a median of 4277 iterations while the solftplus and algebraic penalty functions converged in approximately 844 iterations, a five fold reduction in computational effort.

\begin{figure}
  \includegraphics[width=\linewidth]{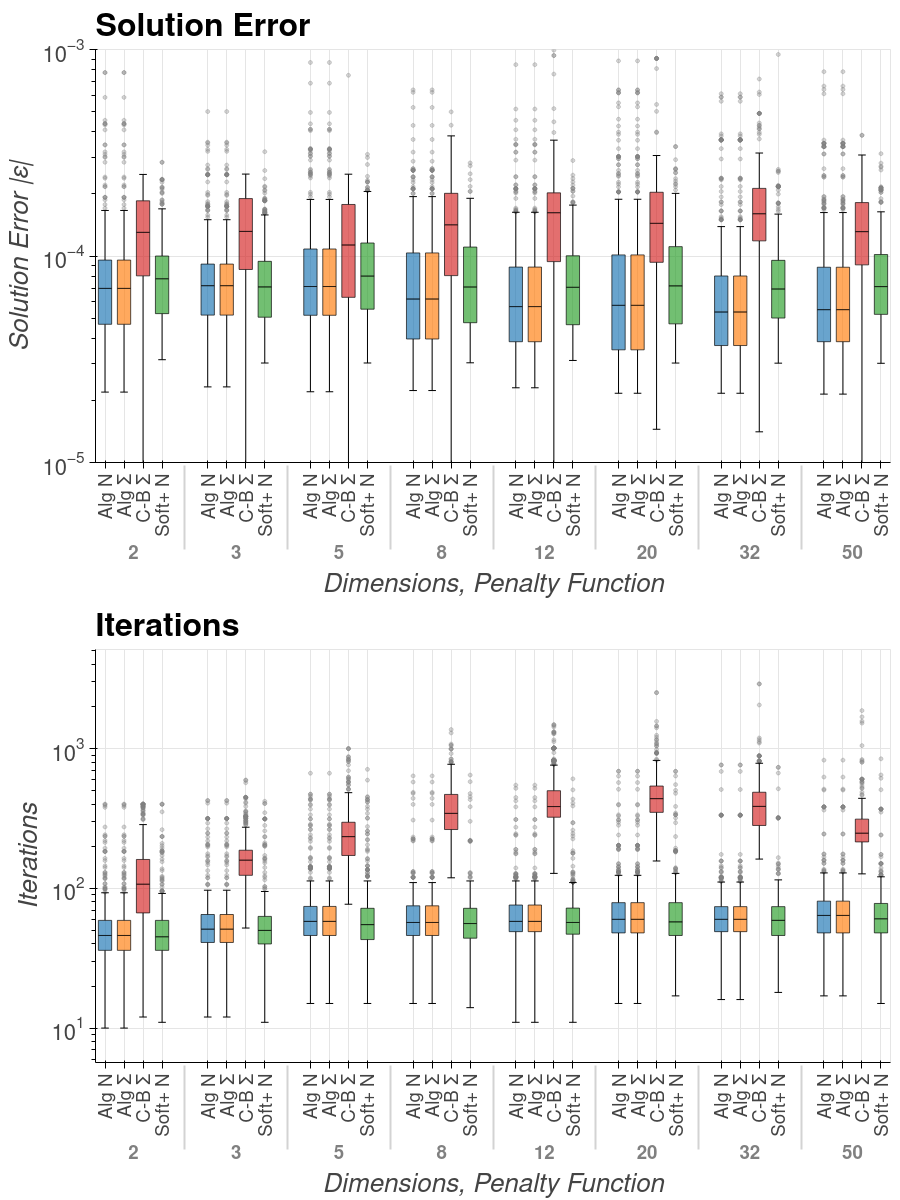}
  \caption{Results - Hyperspheres}
  \label{fig:exp_dims_sphere}
\end{figure}

Table \ref{tab:exp} and Figure \ref{fig:exp_dims_sphere} present the results of the second set of experiments on the hypersphere problem. Solution error is shown to be independent of the number of dimensions for all penalty functions with algebraic and softplus penalty functions returning approximately half as much error than the Courant-Beltrami penalty function. Against a single constraint, the optimization algorithm required between 2 and 7 times as many iterations to converge when using the Courant-Beltrami penalty functions.

A slight advantage in solution error is seen in both experiments when using algebraic penalty functions over the softplus functions. This validates the relationship seen in Figure \ref{fig:norm_err}. In this experiment, $\sigma =$ 15, $\alpha =$ 3e-5, and  $\nabla O(x)$ was randomly selected between 1e-2 and 5. As a result, $|\nabla O(x) / \sigma|$ generally falls in the region where the algebraic function outperforms the softplus function.

\begin{table*}
\caption{Measured Median Solution Error and Number of Iterations}\label{tab:exp}
\center
\resizebox{\linewidth}{!}{
\begin{tabular}{ c | c | c c | c c | c c | c c }
  \multirow{2}{*}{\textbf{Problem}} &\multirow{2}{*}{\textbf{Dimensions}} & \multicolumn{2}{c |}{\textbf{Algebraic Norm}} & \multicolumn{2}{c |}{\textbf{Algebraic Sum}} & \multicolumn{2}{c |}{\textbf{Courant-Beltrami}}& \multicolumn{2}{c}{\textbf{Softplus Norm}}  \\
  && Error & Iterations & Error & Iterations & Error & Iterations & Error & Iterations \\
  \hline
   \multirow{8}{*}{\rotatebox[origin=c]{90}{\textbf{Hyperplanes}}} & 2 & 1.22E-04 & 19 & 1.73E-04 & 18 & 1.89E-04 & 11 & 1.32E-04 & 18 \\
  & 3 & 1.74E-04 & 35 & 2.72E-04 & 33 & 1.60E-04 & 27 & 1.80E-04 & 33 \\
  & 5 & 2.45E-04 & 69 & 5.05E-04 & 79 & 2.39E-04 & 104 & 2.43E-04 & 65 \\
  & 8 & 3.45E-04 & 122.5 & 8.25E-04 & 137 & 4.17E-04 & 314 & 3.81E-04 & 124 \\
  & 12 & 3.95E-04 & 194 & 1.01E-03 & 194.5 & 5.39E-04 & 660 & 4.11E-04 & 187 \\
  & 20 & 6.20E-04 & 378 & 2.02E-03 & 453 & 6.14E-04 & 1602.5 & 6.37E-04 & 364 \\
  & 32 & 7.17E-04 & 557.5 & 2.60E-03 & 686 & 7.77E-04 & 2913.5 & 7.18E-04 & 521 \\
  & 50 & 8.51E-04 & 863 & 3.69E-03 & 1060 & 9.81E-04 & 4277 & 9.07E-04 & 825 \\
  \hline
  \multirow{8}{*}{\rotatebox[origin=c]{90}{\textbf{Hypersphere}}} & 2 & 6.97E-05 & 46 & 6.97E-05 & 46 & 1.30E-04 & 107 & 7.76E-05 & 45 \\
  & 3 & 7.18E-05 & 51 & 7.18E-05 & 51 & 1.32E-04 & 159 & 7.08E-05 & 50 \\
  & 5 & 7.12E-05 & 58 & 7.12E-05 & 58 & 1.13E-04 & 234.5 & 7.99E-05 & 55 \\
  & 8 & 6.19E-05 & 57 & 6.19E-05 & 57 & 1.42E-04 & 344 & 7.08E-05 & 56 \\
  & 12 & 5.69E-05 & 58 & 5.69E-05 & 58 & 1.62E-04 & 385 & 7.05E-05 & 57 \\
  & 20 & 5.77E-05 & 60 & 5.77E-05 & 60 & 1.44E-04 & 438.5 & 7.16E-05 & 57.5 \\
  & 32 & 5.35E-05 & 60 & 5.35E-05 & 60 & 1.60E-04 & 386 & 6.92E-05 & 59 \\
  & 50 & 5.50E-05 & 64 & 5.50E-05 & 64 & 1.31E-04 & 248 & 7.11E-05 & 60.5 \\
\end{tabular}}
\end{table*}

\begin{figure}
  \includegraphics[width=\linewidth]{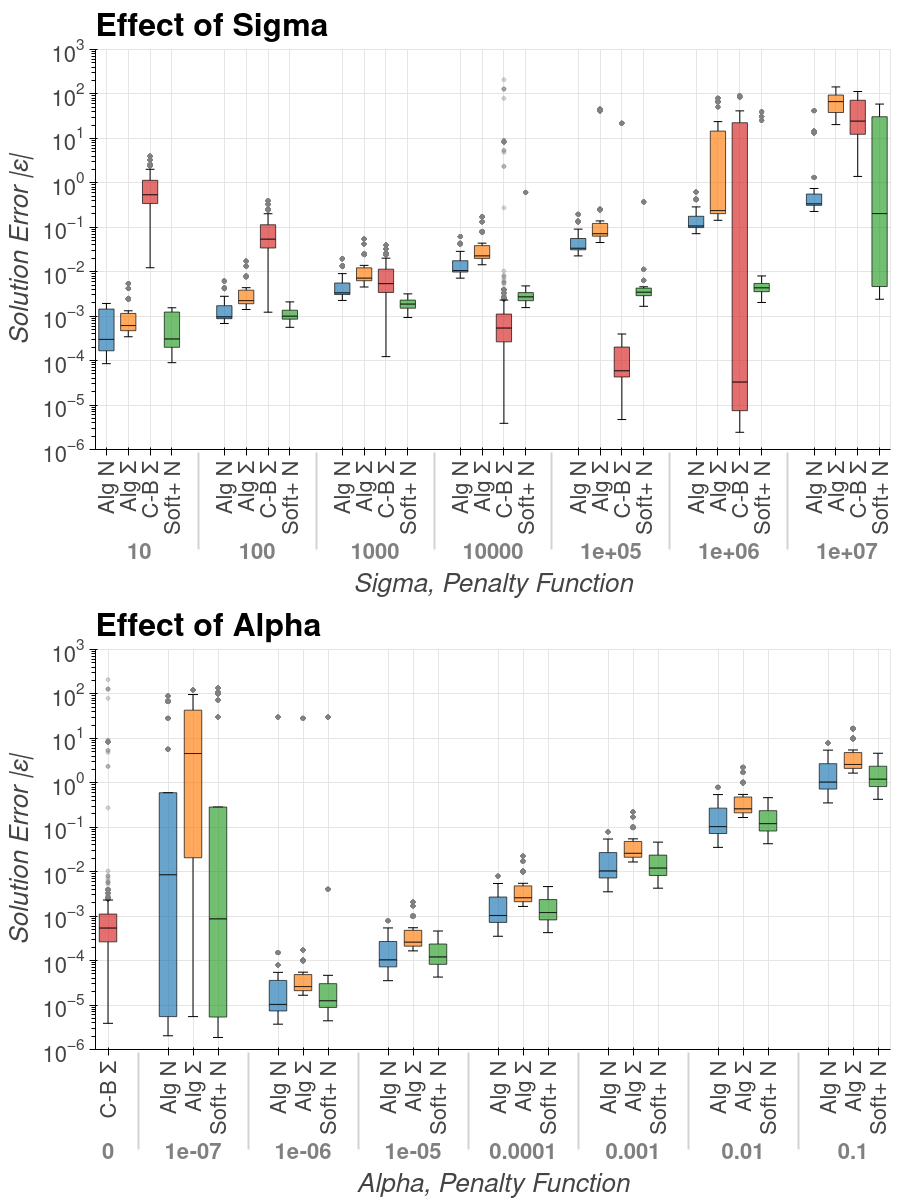}
  \caption{Sensitivity of Error to $\sigma$ and $\alpha$}
  \label{fig:experiment_sigma_alpha}
\end{figure}

The sensitivity of the four penalty function configurations to $\sigma$ and $\alpha$ in the sheared hyperplane experiment is found in Figure \ref{fig:experiment_sigma_alpha}. The solution error measured when using the Courant-Beltrami penalty functions is inversely proportional to $\sigma$ as predicted by Equation \ref{eq:err_cb}. When $\sigma$ approaches 1e6, the inverse of the interval used to approximate the jacobian, the optimization algorthim cannot converge and solutions with large errors are returned.

Solution error of the algebraic and softplus penalty functions increases with $\sigma$, as predicted by Equations \ref{eq:err_alg<>} to \ref{eq:err_log=}. Both functions are relatively insensitive to the selection of $\sigma$, with the softplus penalty function proportional to $log_2(1/\sigma)$, and the algebraic penalty function proportional to $\sqrt{1/\sigma}$. These penalty functions begin to lose effectiveness below $\sigma =$ 10, as the slope of the response to large errors approaches the gradient of the objective function. These functions are also more forgiving with respect to the interval used to approximate the jacobian, and can converge to a reasonable solution for $\sigma >$ 1e6.

The solution error measured when using the algebraic and softplus penalty functions is linearly proportional to the selection of $\alpha$. Converged solutions are possible for values as small as 1e-6, the gradient approximation interval. Reducing interval can allow $\alpha$ to be reduced further.

These experimental results show clear advantages of using algebraic and softplus penalty functions over the Courant-Beltrami penalty function. On the sheared hyperplane problem with 50 dimensions, applying the softplus and algebraic penalty functions reduced the number of iterations required to arrive at a converged solution by a factor of 5. On the hypersphere problem, the new penalty functions converged in one fourth as many iterations to solutions with significantly lower error.

The softplus penalty functions may be preferred for problems where the gradient is not know with confidence, due the low sensitivity to scaling parameter $\sigma$ seen in Figure \ref{fig:experiment_sigma_alpha}. Where execution speed is critical, the algebraic penalty functions are preferable as they incur lower overhead.

\section{Conclusions}\label{conclusions}

Constrained optimization remains a highly relevant field of research with many practical applications. As such, the improved penalty functions proposed in this work could find wide application.

The softplus and algebraic penalty functions developed here were shown to yield solutions as accurate as established methods with one fifth of the computational effort. Furthermore, the new penalty functions where shown to be less sensitive to objective function gradient and scaling, allowing greater control over solution accuracy.

The potential of these penalty functions was only briefly investigated in this work, and further research towards an algorithm setting the scaling parameter $\sigma$ to twice the objective function gradient and towards application to reinforcement learning problems may be worth pursuing. 

\bibliographystyle{plain}
\bibliography{references}

\begin{thebibliography}{1}

\bibitem{benigo}
{\em Incorporating Second-Order Functional Knowledge for Better Option
  Pricing}, Proceedings of the 13th International Conference on Neural
  Information Processing Systems (NIPS'00). MIT Press, 2000.

\bibitem{practical_opt}
Roger Fletcher.
\newblock {\em Practical Methods of Optimization (2nd ed.)}.
\newblock John Wiley and Sons, New York, NY, 1987.

\bibitem{nedich2020convergence}
Angelia Nedich and Tatiana Tatarenko.
\newblock Convergence rate of a penalty method for strongly convex problems
  with linear constraints, 2020.

\bibitem{evolutionary_comp}
Alice~E. Smith and David~W. Coit.
\newblock {\em Constraint-Handling Techniques - Penalty Functions}, chapter
  C5.2.
\newblock Institute of Physics Publishing and Oxford University Press, Bristol,
  U.K., 1997.

\bibitem{tatarenko2020smooth}
Tatiana Tatarenko and Angelia Nedich.
\newblock A smooth inexact penalty reformulation of convex problems with linear
  constraints, 2020.

\end{thebibliography}

\end{document}